\documentclass[11pt]{article}
\usepackage{amsmath,amssymb,amsfonts,amsthm,graphicx}
\usepackage{hyperref}
\newtheorem{theorem}{Theorem}[section]
\newtheorem{proposition}[theorem]{Proposition}
\newtheorem{corollary}[theorem]{Corollary}
\newtheorem{lemma}[theorem]{Lemma}

\newtheorem{remark}[theorem]{Remark}
\newtheorem{definition}[theorem]{Definition}
\numberwithin{equation}{section}
\title{Metric properties of $N$-continued fractions}
\author{
    Dan Lascu\footnote{e-mail: lascudan@gmail.com.}, \nonumber \\
        }
\begin{document}
\maketitle
\thispagestyle{empty}
\begin{abstract}
A generalization of the regular continued fractions was given by
Burger et al. in 2008 \cite{Burger-2008}. In this paper we give metric properties of this expansion. For the transformation which generates this expansion, its invariant measure and Perron-Frobenius operator are investigated.
\end{abstract}

\section{Introduction}
The modern history of continued fractions started with Gauss who found a natural invariant measure of the so-called \textit{regular continued fraction (or Gauss) transformation}, i.e., $T : [0, 1] \rightarrow [0, 1]$, $T(x) = {1}/{x} - \lfloor{1}/{x}\rfloor$, $x \neq 0$, and $T(0) = 0$.
Here $\left\lfloor \cdot \right\rfloor$ denotes the floor (or entire) function.
Let $G$ be this measure which is called \emph{Gauss measure}.
The Gauss measure of an interval $A \in \mathcal{B}_{[0,1]}$ is $G(A) = (1/\log2) \int_{A} 1/(x+1) dx$, where $\mathcal{B}_{[0, 1]}$ denotes the $\sigma$-algebra of all Borel subsets of $[0, 1]$.
This measure is $T$-invariant in the sense that
$G(T^{-1}(A)) = G (A)$ for any $A \in {\mathcal{B}}_I$.

By the very definition, any irrational $0<x<1$ can be written
as the infinite regular continued fraction
\begin{equation}
x = \displaystyle \frac{1}{a_1+\displaystyle \frac{1}{a_2+\displaystyle \frac{1}{a_3+ \ddots}}} :=[a_1, a_2, a_3, \ldots], \label{1.1}
\end{equation}
where $a_n \in \mathbb{N}_+ : = \left\{1, 2, 3, \ldots\right\}$
\cite{Hincin-1964}.
Such $a_n$'s are called {\it incomplete quotients (or continued fraction digits)} of $x$ and they are given by the formulas $a_1(x)=\lfloor1/x\rfloor$ and $a_{n+1}(x)=a_1(T^n(x))$, where $T^n$ denotes the $n$th iterate of $T$.

Thus, the continued fraction representation conjugates the Gauss transformation and the shift on the space of infinite integer-valued sequence $(a_n)_{n \in \mathbb{N}_+}$.

Other famous probabilists like Paul L\'evy and Wolfgang Doeblin also contributed to what is nowadays called the ``metric theory of continued fractions".

The first problem in the metric theory of continued fractions was Gauss' famous 1812 problem \cite{Brezinski}. In a letter dated 1812, Gauss asked Laplace how fast $\lambda(T^{-n}([0, x]))$ converges to the invariant measure $G([0, x])$, where $\lambda$ denotes the Lebesgue measure on $[0, 1]$. Gauss' question was answered independently in 1928 by Kuzmin \cite{Kuzmin-1928}, and in 1929 by Paul L\'evy \cite{Levy-1929}.

Apart from regular continued fractions, there are many other continued fraction expansions: Engel continued fractions, Rosen expansions, the nearest integer continued fraction, the grotesque continued fractions, $f$-expansions etc. For most of these expansions has been proved Gauss-Kuzmin-L\'evy theorem \cite{IosifescuKalpazidou-2009, IS-2006, L-2013, Schweiger, Sebe-2000, Sebe-2001, Sebe-2002, SebeLascu-2014}

The purpose of this paper is to show and prove some metric properties of $N$-continued fraction expansions introduced by Burger et al. \cite{Burger-2008}.

In Section \ref{section2}, we present the current framework. We show a Legendre-type result and the Brod\'en-Borel-L\'evy formula by using the probability structure of $(a_n)_{n \in \mathbb{N}_+}$ under the Lebesgue measure.
In Section \ref{subsection2.5}, we find the invariant measure $G_N$ of $T_N$ the transformation which generate the $N$-continued fraction expansions.
In Section \ref{section3}, we consider the so-called natural extension of
$([0, 1],{\cal B}_{[0, 1]},G_N, T_N)$ \cite{Schweiger}.
In Section \ref{section4}, we derive its Perron-Frobenius operator
under different probability measures on $([0, 1],{\mathcal{B}}_{[0, 1]})$.
Especially, we derive the asymptotic behavior
for the Perron-Frobenius operator of $([0, 1],{\cal B}_{[0, 1]},G_N, T_N)$.

\section{$N$-continued fraction expansions} \label{section2}
In this paper, we consider a generalization of the Gauss transformation.

\subsection{$N$-continued fraction expansions as dynamical system} \label{subsection2.1}
Fix an integer $N \geq 1$. In \cite{Burger-2008}, Burger et al. proved that any irrational  $0<x<1$ can be written in the form
\begin{equation}
x = \displaystyle \frac{N}{a_1+\displaystyle \frac{N}{a_2+\displaystyle \frac{N}{a_3+ \ddots}}} :=[a_1, a_2, a_3, \ldots]_N, \label{2.1}
\end{equation}
where $a_n$'s are non-negative integers. We will call (\ref{2.1}) the \emph{$N$-continued fraction expansion of $x$}.

This continued fraction is treated as the following dynamical systems.
\begin{definition} \label{def2.1}
Fix an integer $N \geq 1$.
\begin{enumerate}
\item[(i)]
The measure-theoretical dynamical system $(I,{\mathcal B}_{I},T_N)$ is defined as follows:
$I:=[0,1]$,
$\mathcal{B}_I$ denotes the $\sigma$-algebra of all Borel subsets of $I$,
and $T_N$ is the transformation
\begin{equation}
T_N:I \to I; \quad
T_{N}(x):=
\left\{
\begin{array}{ll}
{\displaystyle \frac{N}{x}- \left\lfloor\frac{N}{x}\right\rfloor}&
{ \mbox{if } x \in I,}\\
\\
0& \mbox{if }x=0.
\end{array}
\right. \label{2.2}
\end{equation}
%

\item[(ii)]
In addition to (i),
we write $(I,{\cal B}_{I},G_N,T_N)$ as
$(I,{\cal B}_{I}, T_N)$ with
the following probability measure $G_N$ on $(I,{\cal B}_{I})$:
\begin{equation}
G_N (A) :=
\frac{1}{\log \frac{N+1}{N}} \int_{A} \frac{dx}{x+N},
\quad A \in {\mathcal{B}}_I. \label{2.3}
\end{equation}
\end{enumerate}
\end{definition}

Define the {\it quantized index map} $\eta: I \to {\mathbb N}:={\mathbb N}_+\cup \{0\}$ by
\begin{equation}
\eta(x) := \left\{\begin{array}{ll}
\displaystyle \left\lfloor \frac{N}{x} \right\rfloor & \hbox{if }  x \neq 0, \\
\\
\infty & \hbox{if }  x = 0.
\end{array} \right. \label{2.4}
\end{equation}
By using $T_N$ and $\eta$,
the sequence $(a_{n})$ in (\ref{2.1}) is obtained as follows:
\begin{equation}
a_n= \eta\left(T_N^{n-1}(x)\right), \quad n \geq 1 \label{2.5}
\end{equation}
with $T_N^0 (x) = x$. Since $x \in (0, 1)$ we have that $a_n \geq N$ for any $n \geq 1$.

In this way, $T_N$ gives the algorithm of $N$-continued fraction expansion which is an obvious generalization of the regular continued fraction.

\begin{proposition} \label{prop3.1}
Let $(I,{\cal B}_{I},G_N,T_N)$ be as in Definition \ref{def2.1}(ii).
\begin{enumerate}
\item[(i)]
$(I,{\cal B}_{I},G_N,T_N)$ is ergodic.
\item[(ii)]
The measure $G_N$ is invariant under $T_N$.
\end{enumerate}
\end{proposition}
\noindent \textbf{Proof.}
See \cite{DKW-2013} and Section \ref{subsection2.5}
\hfill $\Box$\\

\noindent
By Proposition \ref{prop3.1}(ii), $(I,{\cal B}_{I},G_N, T_N)$ is a ``dynamical system" in the sense of Definition 3.1.3 in \cite{BG-1997}.

\subsection{Some elementary properties of $N$-continued fractions} \label{subsection2.2}
Roughly speaking, the metrical theory of continued fraction expansions is
the asymptotic analysis of incomplete quotients $(a_n)_{n \in \mathbb{N}_+}$ and related sequences \cite{IK-2002}.
First, note that in the rational case, the continued fraction expansion (\ref{2.1}) is finite, unlike the irrational case, when we have an infinite number of digits.
In \cite{Wekken-2011}, Van der Wekken showed the convergence of the expansion.
For $x \in I \setminus {\mathbb{Q}}$, define the {\it $n$-th order convergent} $[a_1, a_2, \ldots, a_n]_N$ of $x$ by truncating the expansion on the right-hand side of (\ref{2.1}), that is,
\begin{equation}
[a_1, a_2, \ldots, a_n]_N \to x, \quad n \to \infty.  \label{2.6}
\end{equation}
To this end, for $n \in \mathbb{N}_+$, define integer-valued functions $p_n(x)$ and $q_n(x)$ by
\begin{eqnarray}
p_n(x) &:=& a_n p_{n-1} + N p_{n-2}, \quad n \geq 2 \label{2.7} \\
q_n(x) &:=& a_n q_{n-1} + N q_{n-2}, \quad n \geq 1 \label{2.8}
\end{eqnarray}
with $p_0(x):=0$, $q_0(x) := 1$, $p_{-1}(x):=1$, $q_{-1}(x):=0$, $p_1(x):=N$, $q_1(x):=a_1$.
By induction, we have
\begin{equation}
p_{n-1}(x)q_n(x) - p_n(x)q_{n-1}(x) = (-N)^n, \quad n \in \mathbb{N}. \label{2.9}
\end{equation}
By using (\ref{2.7}) and (\ref{2.8}), we can verify that
\begin{equation}
x = \frac{p_n(x) + T^n_N(x)p_{n-1}(x)}
{q_n(x) + T^n_N(x)q_{n-1}(x)}, \quad
n \geq 1. \label{2.10}
\end{equation}
By taking $T^n_N(x)=0$ in (\ref{2.10}), we obtain $ [a_1, a_2, \ldots, a_n]_N = p_n(x)/q_n(x)$. From this and by using (\ref{2.9}) and (\ref{2.10}), we obtain
\begin{equation}
\left| x - \frac{p_n(x)}{q_n(x)} \right| =
\frac{N^n \cdot T^n_N(x)}{q_n(x)\left(T^n_N(x)q_{n-1}(x)+q_n(x)\right)}, \quad n \geq 1. \label{2.11}
\end{equation}
Now, since $T^n_N(x)<1$ and $\left|T^n_N(x)\displaystyle \frac{q_{n-1}(x)}{q_n(x)} + 1 \right| \geq 1$, we have
\begin{equation}
\left| x - \frac{p_n(x)}{q_n(x)} \right| < \frac{N^n}{q^2_n(x)}, \quad n \geq 1. \label{2.12}
\end{equation}
In order to prove (\ref{2.6}), it is sufficient to show the following inequality:
\begin{equation}
\left| x - \frac{p_n(x)}{q_n(x)} \right| \leq \frac{1}{N^n}, \quad n \geq 1. \label{2.13}
\end{equation}
From (\ref{2.8}), we have that $q_n(x) > N q_{n-1}(x)$ and because $q_0=1$ we have $q_n(x) > N^n$. Finally, (\ref{2.13}) follows from (\ref{2.12}).

\subsection{Diophantine approximation} \label{subsection2.3}
Diophantine approximation
deals with the approximation of real numbers by rational numbers \cite{Hincin-1964}.
We approximate $x \in I \setminus {\mathbb{Q}}$ by incomplete quotients in (\ref{2.7}) and (\ref{2.8}).

For $x \in I \setminus {\mathbb{Q}}$, let $a_n$ be as in (\ref{2.5}).
For any $n \in \mathbb{N}_+$ and $i^{(n)}=(i_1, \ldots, i_n) \in \mathbb{N}^n$,
define the {\it fundamental interval associated with} $i^{(n)}$ by
\begin{equation}
I_N \left(i^{(n)}\right) = \{x \in I \setminus {\mathbb{Q}}:  a_k(x) = i_k
\mbox{ for } k=1,\ldots,n \} \label{2.14}
\end{equation}
where we write $I_N (i^{(0)}) = I \setminus {\mathbb{Q}}$.
Remark that $I_N (i^{(n)})$ is not connected by definition.
For example, we have
\begin{equation}
I_N(i) \,= \,\{x \in I \setminus {\mathbb{Q}}: a_1 = i \}
\,=\,\left( I \setminus {\mathbb{Q}}\right) \cap \left(\frac{N}{i+1}, \frac{N}{i}\right)
\quad
\mbox{for any }i \in \mathbb{N}. \label{2.15}
\end{equation}

\begin{lemma} \label{lema2.3}
Let $\lambda$ denote the Lebesgue measure. Then
\begin{equation}
\lambda \left(I_N (i^{(n)})\right) = \frac{N^n}{q_n(x)(q_n(x)+q_{n-1}(x))},
\label{2.16}
\end{equation}
where $(q_n)$ is as in (\ref{2.8}).
\end{lemma}
\noindent \textbf{Proof.}
From the definition of $T_N$ and (\ref{2.10}), we have
\begin{equation}
I_N\left(i^{(n)}\right)
\,= \,\left( I \setminus {\mathbb{Q}}\right) \,\cap \,\left(\,u(i^{(n)}), \,v(i^{(n)})\,\right), \label{2.17}
\end{equation}
where both $u\left(i^{(n)}\right)$ and $v\left(i^{(n)}\right)$
are rational numbers defined as
\begin{equation}
u\left(i^{(n)}\right):=\left\{
\begin{array}{lll}
	\displaystyle\frac{p_n(x)+p_{n-1}(x)}{q_n(x)+q_{n-1}(x)}
& \quad \mbox{if $n$ is odd}, \\
	\\
	\displaystyle\frac{p_n(x)}{q_n(x)} & \quad \mbox{if $n$ is even,} \\
\end{array}
\right. \label{2.18}
\end{equation}
and
\begin{equation}
v\left(i^{(n)}\right):=\left\{
\begin{array}{lll}
	\displaystyle\frac{p_n(x)}{q_n(x)} & \quad \mbox{if $n$ is odd,} \\
	\\
	\displaystyle\frac{p_n(x)+p_{n-1}(x)}{q_n(x)+q_{n-1}(x)} & \quad \mbox{if $n$ is even}.
\end{array}
\right. \label{2.19}
\end{equation}
By using (\ref{2.9}), we have (\ref{2.16}).
\hfill $\Box$\\

We now give a Legendre-type result for $N$-continued fraction expansions.
For $x \in I \setminus {\mathbb{Q}}$, we define the \textit{approximation coefficient} $\Theta_N(x,n) $ by
\begin{equation}
\Theta_N(x,n)
:= \frac{q_n^2}{N^n} \left|x-\frac{p_n}{q_n}\right|, \quad n \geq 1 \label{2.20}
\end{equation}
where $p_n / q_n$ is the $n$th continued fraction convergent of $x$ in (\ref{2.1}).

\begin{proposition} \label{prop.2.4}
For $x \in I \setminus {\mathbb{Q}}$ and an irreducible fraction $0<p/q<1$, assume that $p/q$ is written as follows:
\begin{equation} \label{2.21a}
\frac{p}{q} = [i_1, \ldots, i_n]_N
\end{equation}
where $[i_1, \ldots, i_n]_N$ is as in $(\ref{2.6})$, and the length
$n \in \mathbb{N}_+$ of $N$-continued fraction expansion of \mbox{ }$p/q$ is chosen in such a way that it is even if \mbox{ }$p/q < x$ and odd otherwise.
Then
\begin{equation} \label{2.21}
\Theta_N(x,n)
< \frac{q}{q+q_{n-1}} \quad \mbox{if and only if} \quad  \frac{p}{q} \mbox { is the $n$th convergent of } x
\end{equation}
where $\Theta_N(x,n)$ is as in (\ref{2.20}) and the positive integer $q_{n-1}$ is defined as the denominator  of the irreducible fraction representation of the rational number $[i_1, \ldots, i_{n-1}]_N$
with $q_0=1$ for the sequence $i_1, \ldots, i_n$.
\end{proposition}

\noindent \textbf{Proof.}
Fix $x \in I \setminus {\mathbb{Q}}$ and $n \geq 1$. Let $\Theta:=\Theta_N(x,n)$.

\noindent
($\Leftarrow$)
Assume that $p/q$ is the $n$th convergent of $x$.
By (\ref{2.11}) and the definition of $\Theta$,
we have
\begin{equation} \label{2.23}
\Theta = \frac{q^2}{N^n} \left|x - \frac{p}{q}\right| = \frac{T^n_N(x)q}{q+T^n_N(x)q_{n-1}(x)} \leq \frac{q}{q+q_{n-1}}
\end{equation}
where we use $q_{n-1}=q_{n-1}(x)$.

\noindent
($\Rightarrow$)
Conversely,
\begin{equation} \label{2.24}
\mbox{if }\Theta < \displaystyle\frac{q}{q+q_{n-1}}, \quad \mbox{ then } \quad
q\left|x-\frac{p}{q}\right| < \frac{N^n}{q+q_{n-1}}.
\end{equation}
If $n$ is even, then $x > p/q$ and we have
\begin{equation} \label{2.25}
x-\frac{p}{q} <
\frac{N^n}{q(q+q_{n-1})}.
\end{equation}
From these,
\begin{equation} \label{2.26}
\frac{p}{q}< x < \frac{p}{q} + \frac{N^n}{q(q+q_{n-1})} = \frac{p+p_{n-1}}{q+q_{n-1}}
\end{equation}
where $p_{n-1}$ is defined as
$p_{n-1}/q_{n-1}=[i_1, \ldots, i_{n-1}]_N$.
Hence $x \in I_N \left(i^{(n)}\right)$, i.e.,
$p/q = [i_1, \ldots, i_n]_N$ is a convergent of $x$. The case when $n$ is odd is treated similarly.
\hfill $\Box$

\subsection{Brod\'en-Borel-L\'evy formula and its consequences} \label{subsection2.4}
We derive the so-called Brod\'en-Borel-L\'evy formula \cite{IG-2009, IK-2002} for $N$-continued fraction expansion.
For $x \in I$, let $a_n$ and $q_n$ be as in (\ref{2.5}) and (\ref{2.8}), respectively.
We define $(s_n)_{n \geq 0}$ by
\begin{equation}
s_0:=0,\quad
s_n := N \frac{q_n}{q_{n-1}}, \quad n \geq 1. \label{2.27}
\end{equation}
From (\ref{2.8}), $s_n = N/(a_n+s_{n-1})$ for $n \geq 1$. Hence
\begin{equation}
s_n = \frac{N}{\displaystyle a_n+\frac{N}{\displaystyle a_{n-1}+\ddots+ \frac{N}{a_1}}} = [a_n, a_{n-1}, \ldots, a_2, a_1]_N, \label{2.28}
\end{equation}
for $n \geq 1$.
For a probability measure $\mu$, recall that the {\it conditional probability} of an event $A$ assuming that $B$ has occurred, denoted $\mu(A|B)$, is defined as $\mu(A|B) := \mu(A \cap B )/\mu(B)$, $\mu(B)>0$.
\begin{proposition} [Brod\'en-Borel-L\'evy formula] \label{prop2.5}
Let $\lambda$ denote the Lebesgue measure on $I$.
For any $n \in \mathbb{N}_+$, the conditional probability
$\lambda (T^n_N < x |a_1,\ldots, a_n )$ is given as follows:
\begin{equation}
\lambda (T^n_N < x |a_1,\ldots, a_n )
= \frac{(s_n + N)x}{s_n x+N}, \quad x \in I \label{2.29}
\end{equation}
where $s_n$ is as in $(\ref{2.27})$ and $a_1, \ldots, a_n$ are as in $(\ref{2.5})$.
\end{proposition}
\noindent\textbf{Proof.}
By definition, we have
\begin{equation} \label{2.30}
\lambda\left(T^n_N < x |a_1,\ldots, a_n \right) =
\frac{\lambda\left(\left(T^n_N < x\right) \cap I_N(a_1,\ldots, a_n) \right)}{\lambda\left(I_N(a_1,\ldots, a_n)\right)}
\end{equation}
for any $n \in \mathbb{N}_+$ and $x \in I$.
From (\ref{2.10}) and (\ref{2.17}) we have
\begin{eqnarray} \label{2.31}
\lambda\left(\left(T^n_N < x \right) \cap I_N(a_1,\ldots, a_n)\right) &=& \left|\frac{p_n}{q_n} - \frac{p_n+xp_{n-1}}{q_n+xq_{n-1}}\right|  \\
\nonumber \\
&=& \frac{N^nx}{q_n(q_n+xq_{n-1})}. \nonumber
\end{eqnarray}
From this and (\ref{2.16}), we have
\begin{eqnarray} \label{2.32}
\lambda\left(T^n_N < x |a_1, \ldots, a_n \right) &=& \frac{\lambda\left(\left(T^n_N < x\right) \cap I_N(a_1,\ldots, a_n) \right)}{\lambda\left(I_N(a_1,\ldots, a_n)\right)} \nonumber \\
\nonumber \\
&=& \frac{x\left(q_n+q_{n-1}\right)}{q_n+xq_{n-1}} \nonumber \\
\nonumber \\
&=& \frac{(s_n+N)x}{s_n x+N}
\end{eqnarray}
for any $n \in \mathbb{N}_+$ and $x \in I$.
\hfill $\Box$

The Brod\'en-Borel-L\'evy formula allows us to determine the probability structure of incomplete quotients $(a_n)_{n \in \mathbb{N}_+}$ under $\lambda$.
\begin{proposition} \label{prop2.6}
For any $i \geq N$ and $n \in \mathbb{N}_+$, we have
\begin{equation}
\lambda(a_1=i) = \frac{N}{i(i+1)}, \quad
\lambda\left(a_{n+1}=i |a_1,\ldots, a_n \right) = V_{N,i}(s_n) \label{2.33}
\end{equation}
where $(s_n)$ is as in (\ref{2.27}), and
\begin{equation}
V_{N,i}(x) := \frac{x+N}{(x+i)\,(x+i+1)}. \label{2.34}
\end{equation}
\end{proposition}
\noindent\textbf{Proof.}
From (\ref{2.15}), the case $\lambda(a_1=i)$ holds.
For $n \geq N$ and $x \in I \setminus {\mathbb{Q}}$, we have
$ T_N^n(x) = [a_{n+1}, a_{n+2}, \ldots]_N$
where $(a_n)$ is as in (\ref{2.5}).
By using (\ref{2.29}), we have
\begin{eqnarray} \label{2.35}
\lambda(\,a_{n+1}=i \,|\,a_1,\ldots, a_n\, )
& = & \lambda\left( \,T^n_N \in \left(\frac{N}{i+1},\frac{N}{i} \right]\,
| \,a_1,\ldots, a_n \,\right). \nonumber \\
&=& \frac{(s_n+N)\frac{N}{i}}{s_n \frac{N}{i}+N} - \frac{(s_n+N)\frac{N}{i+1}}{s_n \frac{N}{i+1}+N} \nonumber \\
\nonumber\\
&=& V_{N,i}(s_n).
\end{eqnarray}
\hfill $\Box$\\
In (\ref{2.33}),
$\sum_{i=N}^{\infty}\lambda(a_{n+1}=i |a_1,\ldots, a_n )$
must be $1$ because $\lambda$ is a probability measure on $(I,{\cal B}_I)$.
This can be verified from (\ref{2.33}) and (\ref{2.35}) by using the partial fraction decomposition. By the same token, we see that
\begin{equation} \label{2.36}
\sum_{i=N}^{\infty}V_{N,i}(x)=1 \quad \mbox{ for any }x \in I.
\end{equation}

\begin{remark}
Proposition (\ref{prop2.6}) is the starting point of an approach to the metrical
theory of $N$-continued fraction expansions via dependence with complete connections (see
\cite{IG-2009}, Section 5.2) which will the subject of our future paper.
\end{remark}

\begin{corollary}
The sequence $(s_n)_{n \in \mathbb{N}_+}$ with $s_0 = 0$ is a homogeneous $I$-valuated Markov chain on $(I, \mathcal{B}_I, \lambda)$ with the following transition mechanism: from state $s \in I$ the only possible one-step transitions are those to states $N/(s+i)$, $i \geq N$, with corresponding  probabilities $V_{N,i}(s)$, $i \geq N$.
\end{corollary}

\subsection{The invariant measure of $T_N$} \label{subsection2.5}
Let $(I, \mathcal{B}_I)$ be as in Definition \ref{def2.1}(i). In this subsection we will give the explicit form of the invariant probability measure $G_N$ of the transformation $T_N$ in (\ref{2.2}), i.e., $G_N((T_N)^{-1}(A)) = G_N(A)$ for any $A \in {\mathcal{B}}_I$.
From the aspect of metric theory, the digits $a_n$ in (\ref{2.5}) can be viewed as random variables on $(I, \mathcal{B}_I)$ that are defined almost surely with respect to any probability measure on $\mathcal{B}_I$ assigning probability $0$ to the set of rationals in $I$. Such a probability measure is Lebesgue measure $\lambda$, but a more important one in the present context is the invariant probability measure $G_N$ of the transformation $T_N$.

\begin{proposition} \label{prop2.9}
  The invariant probability density $\rho_N$ of the transformation $T_N$ is given by
\begin{equation}\label{2.37}
  \rho_N(x)=\frac{k_N}{x+N}, \quad x \in I
\end{equation}
where $k_N$ is the normalized constant such that the invariant measure $G_N$ is a probability measure. Furthermore, the constant $k_N$ is given in (\ref{2.3}), i.e., $k_N=\left(\log\left(\frac{N+1}{N}\right)\right)^{-1}$.
\end{proposition}
\noindent\textbf{Proof.}
We will give a proof which involves properties of the Perron-Frobenius operator of $T_N$ under $G_N$. Therefore, the proof will be given in Section \ref{section4}.

\section{Natural extension and extended random variables} \label{section3}
Fix an integer $N \geq 1$. In this section, we introduce the natural extension $\overline{T_N}$ of $T_N$ in (\ref{2.2}) and its extended random variables according to Chap. 21.2 of \cite{Schweiger}.

\subsection{Natural extension} \label{subsection3.1}
Let $(I,{\mathcal B}_I,T_N)$ be as in Definition \ref{def2.1}(i).
Define $(u_{N,i})_{i \geq N}$ by
\begin{equation}
u_{N,i}:I \rightarrow I; \quad
u_{N,i}(x) := \frac{N}{x+i}, \quad x \in I. \label{3.1}
\end{equation}
\noindent
For each $i \geq N$, $u_{N,i}$ is a right inverse of $T_N$, that is,
\begin{equation} \label{3.2}
(T_N \circ u_{N,i})(x) = x, \quad \mbox{for any } x \in I.
\end{equation}
Furthermore, if $\eta(x)=i$, then $(u_{N,i} \circ T_N)(x)=x$ where $\eta$ is as in (\ref{2.4}).
\begin{definition} \label{def3.1}
The \textit{natural extension} $(I^2, {\mathcal B}_{I^2},\overline{T_N})$ of $(I,{\mathcal B}_I,T_N)$ is the transformation $\overline{T_N}$
of the square space
$(I^2,{\mathcal B}_I^2):=(I,{\cal B}_I) \times (I,{\cal B}_I)$
defined as follows  \cite{Nakada-1981}:
\begin{equation} \label{3.3}
\overline{T_N}: I^2 \rightarrow I^2; \quad
\overline{T_N}(x,y) := \left( T_N(x), \,u_{N,\eta(x)}(y) \right), \quad (x, y) \in I^2.
\end{equation}
\end{definition}
\noindent
From (\ref{3.2}), we see that $\overline{T_N}$ is bijective on $I^2$ with the inverse
\begin{equation} \label{3.4}
(\overline{T_N})^{-1}(x, y)
= (u_{N,\eta(y)}(x), \,
T_N(y)), \quad (x, y) \in I^2.
\end{equation}
Iterations of (\ref{3.3}) and (\ref{3.4}) are given as follows for each $n \geq 2$:
\begin{eqnarray}
\left(\overline{T_N}\right)^n(x, y) &=&
\left(\,T^n_N(x), \,[x_n, x_{n-1}, \ldots, x_2(x),\, x_1+ y ]_N \,\right), \label{3.5} \\
\nonumber
\\
\left(\overline{T_N}\right)^{-n}(x, y) &=&
\left(\,[y_n, y_{n-1}, \ldots, y_2, \,y_1+x ]_N,\, T^{n}_N(y) \,\right)
\label{3.6}
\end{eqnarray}
where
$x_i:=\eta\left(T_N^{i-1}(x)\right)$ and
$y_i:=\eta\left(T_N^{i-1}(y)\right)$ for $i=1,\ldots,n$.

For $G_N$ in (\ref{2.3}), Dajani et al. \cite{DKW-2013} define its {\it extended measure} $\overline{G_N}$ on $(I^2, {\mathcal{B}}^2_I)$ as
\begin{equation} \label{3.7}
\overline{G_N}(B) := k_N \int\!\!\!\int_{B}
\frac{dxdy}{(xy+N)^2}, \quad B \in {\mathcal{B}}^2_I.
\end{equation}
Then
$\overline{G_N}(A \times I) = \overline{G_N}(I \times A) = G_N(A)$ for any $A \in {\mathcal{B}}_I$.
The measure $\overline{G_N}$ is preserved by $\overline{T_N}$ \cite{DKW-2013}.

\subsection{Extended random variables} \label{subsection3.2}
Define the projection $E:I^2 \rightarrow I$
by $E(x,y):=x$.
With respect to $\overline{T_N}$ in (\ref{3.3}),
define {\it extended incomplete quotients} $\overline{a}_l(x,y)$,
$l \in \mathbb{Z}$ at $(x, y) \in I^2$ by
\begin{equation} \label{3.8}
\overline{a}_{l}(x, y) := (\eta \circ E)\left(\,(\overline{T_N})^{l-1} (x, y) \,\right),
\quad l \in \mathbb{Z}.
\end{equation}
\begin{remark} \label{rem3.2}
{\rm
\begin{enumerate}
\item[(i)]
Remark that $\overline{a}_{l}(x, y)$ in (\ref{3.8})
is also well-defined for $l \leq 0$ because $\overline{T_N}$ is invertible.
By (\ref{3.5}) and (\ref{3.6}), we have
\begin{equation} \label{3.9}
\overline{a}_n(x, y) = x_n, \quad
\overline{a}_0(x, y) = y_1, \quad
\overline{a}_{-n}(x, y) = y_{n+1},
\quad n \in \mathbb{N}_+, \ (x, y) \in I^2
\end{equation}
where we use notations in  (\ref{3.5}) and (\ref{3.6}).
\item[(ii)]
Since the measure $\overline{G_N}$ is preserved by $\overline{T_N}$, the doubly infinite sequence $(\overline{a}_l(x,y))_{l \in \mathbb{Z}}$
is strictly stationary (i.e., its distribution is invariant under a shift of the indices) under $\overline{G_N}$.
\end{enumerate}
}
\end{remark}
\begin{theorem} \label{th3.3}
Fix $(x,y) \in I^2$ and let $\overline{a}_{l}:=\overline{a}_l(x,y)$ for $l \in {\mathbb Z}$.
Define $a:= [\overline{a}_0, \overline{a}_{-1}, \ldots]_N$.
Then the following holds for any $x \in I$:
\begin{equation} \label{3.10}
\overline{G_N} ( [0, x] \times I \,|
\,\overline{a}_0, \overline{a}_{-1}, \ldots )
= \frac{(N+a)x}{ax + N} \quad \overline{G_N} \mbox{-}\mathrm{a.s.}
\end{equation}
\end{theorem}
\noindent \textbf{Proof.} Recall fundamental interval in (\ref{2.14}).
Let $I_{N,n}$ denote the fundamental interval $I_N(\overline{a}_0, \overline{a}_{-1}, \ldots, \overline{a}_{-n})$ for $n \in \mathbb{N}$. We have
\begin{equation} \label{3.11}
\overline{G_N} ( [0, x] \times I \left. \right| \overline{a}_0, \overline{a}_{-1}, \ldots ) =
\lim_{n \rightarrow \infty} \overline{G_N} ( [0, x] \times I \left. \right| \overline{a}_0, \ldots, \overline{a}_{-n} ) \quad \overline{G_N} \mbox{-a.s.}
\end{equation}
and
\begin{eqnarray} \label{3.12}
\overline{G_N} ( [0, x] \times I \left. \right| \overline{a}_0, \ldots, \overline{a}_{-n} )
&=&
\displaystyle{\frac{\overline{G_N} ([0, x] \times I_{N,n})}{\overline{G_N} (I \times I_{N,n})}} \nonumber \\
\nonumber \\
 &=& \displaystyle{\frac{k_m}{G_N(I_{N,n})}\int_{I_{N,n}}dy\displaystyle\int^x_0{\frac{Ndu}{(yu+N)^2}}} \nonumber \\
 \nonumber \\
 &=& \displaystyle{\frac{1}{G_N(I_{N,n})} \int_{I_{N,n}} \frac{x(y+N)}{xy+N}\, G_N(dy)} \nonumber \\
 \nonumber \\
 &=& \displaystyle{\frac{x(y_n+N)}{xy_n+N}}
\end{eqnarray}
for some $y_n \in I_{N,n}$ where $k_m$ is as in Proposition (\ref{prop2.9}). Since
\begin{equation} \label{3.13}
\lim_{n \rightarrow \infty} y_n =[\overline{a}_0, \overline{a}_{-1}, \ldots]_N = a, \end{equation}
the proof is completed.
\hfill $\Box$

The stochastic property of $(\overline{a}_l)_{l \in \mathbb{Z}}$ under $\overline{G_N}$ is given as follows.
\begin{corollary} \label{cor3.4}
For any $i \in \mathbb{N}$, we have
\begin{equation} \label{3.14}
\overline{G_N} (\left.\overline{a}_1 = i\right| \overline{a}_0, \overline{a}_{-1}, \ldots) = V_{N,i}(a) \quad \overline{G_N} \mbox{-}\mathrm{a.s.}
\end{equation}
where $a = [\overline{a}_0, \overline{a}_{-1}, \ldots]_N$ and $V_{N,i}$ is as in $(\ref{2.34})$.
\end{corollary}
\noindent \textbf{Proof.} Let $I_{N,n}$ be as in the proof of Theorem \ref{th3.3}.
We have
\begin{equation} \label{3.15}
\overline{G_N} (\left.\overline{a}_1 = i\,\right| \,
\overline{a}_0, \overline{a}_{-1}, \ldots) = \lim_{n \rightarrow \infty}
\overline{G_N} (\left.\overline{a}_1 = i\,\right| \,I_{N,n}).
\end{equation}
Now
\begin{eqnarray} \label{3.16}
\overline{G_N} \left( \left.\left( \frac{N}{i+1}, \frac{N}{i} \right] \times
 [0, 1)\right| I_{N,n}\right) &=&
\frac{\overline{G_N} \left(\left( \frac{N}{i+1}, \frac{N}{i} \right] \times I_{N,n}\right)}
{\overline{G_N} (I \times I_{N,n})} \nonumber \\
\nonumber \\
& = & \frac{1}{G_N (I_{N,n})} \int_{I_{N,n}} V_{N,i}(y)\, G_N(dy) \nonumber \\
\nonumber \\
& = & V_{N,i}(y_n)
\end{eqnarray}
for some $y_n \in I_{N,n}$. From (\ref{3.13}), the proof is completed.
\hfill $\Box$
\begin{remark} \label{rem3.5}
\rm The strict stationarity of $\left(\overline{a}_l\right)_{l \in \mathbb{Z}}$, under $\overline{G_N}$ implies that
\begin{equation} \label{3.17}
\overline{G_N}(\left.\overline{a}_{l+1} = i\, \right|\, \overline{a}_l,
\overline{a}_{l-1}, \ldots)
= V_{N,i}(a) \quad \overline{G_N} \mbox{-}\mathrm{a.s.}
\end{equation}
for any $i \in \mathbb{N}$ and $l \in \mathbb{Z}$, where
$a = [\overline{a}_l, \overline{a}_{l-1}, \ldots]_N$.
The last equation emphasizes that $\left(\overline{a}_l\right)_{l \in \mathbb{Z}}$ is a chain of infinite order in the theory of dependence with complete connections (see \rm{\cite{IG-2009}}, Section 5.5).
\end{remark}

\section{Perron-Frobenius operators} \label{section4}

Let $(I,{\mathcal B}_{I},G_N, T_N)$ be as in Definition \ref{def2.1}(ii).
In this section, we derive its Perron-Frobenius operator.

Let $\mu$ be a probability measure on $(I, {\mathcal{B}}_I)$
such that $\mu((T_N)^{-1}(A)) = 0$ whenever $\mu(A) = 0$ for
$A \in {\mathcal{B}}_{I}$.
For example, this condition is satisfied if $T_N$ is $\mu$-preserving, that is, $\mu (T_N)^{-1} = \mu$.
Let %
$
L^1(I, \mu):=\{f: I \rightarrow \mathbb{C} : \int_{I} |f |d\mu < \infty \}.
$
The {\it Perron-Frobenius operator}
of $(I,{\mathcal B}_I,\mu,G_N)$
is defined as the bounded linear operator $U$
on the Banach space $L^1(I,\mu)$ such that the following holds:
\begin{equation}
\int_{A}Uf \,d\mu = \int_{(T_{N})^{-1}(A)}f\, d\mu \quad
\mbox{ for all }
A \in {\mathcal{B}}_{I},\, f \in L^1(I,\mu). \label{4.1}
\end{equation}
About more details, see \cite{BG-1997, IK-2002} or Appendix A in \cite{L-2013}.
\begin{proposition} \label{prop4.1}
Let $(I,{\mathcal B}_{I},G_{N}, T_{N})$ be as in Definition \ref{def2.1}(ii), and let $U$ denote its Perron-Frobenius operator. Then the following holds:
\begin{enumerate}
\item[(i)]
The following equation holds:
\begin{equation}
Uf(x) = \sum_{i \geq N}V_{N,i}(x)\,f\left(\frac{N}{x+i}\right), \quad
f \in L^1(I,G_{N}), \label{4.2}
\end{equation}
where $V_{N,i}$ is as in $(\ref{2.34})$.
\item[(ii)]
Let $\mu$ be a probability measure on $(I,{\mathcal{B}}_I)$
such that $\mu$ is absolutely continuous with respect to
the Lebesgue measure $\lambda$
and let $h := d\mu / d \lambda$  a.e. in $I$.
Then the following holds:
\begin{enumerate}
	\item[(a)] Let $S$ denote the Perron-Frobenius operator of \, $T_N$ under $\mu$.
Then the following holds a.e. in $I$:
	\begin{eqnarray}
S f(x) &=& \frac{N}{h(x)} \sum_{i \geq N} \frac{h\left(\frac{N}{x+i}\right)}
{(x+i)^2} f\left(\frac{N}{x+i}\right) \label{4.3} \\
             &=& \frac{U \hat{f}(x)}{(x+N)h(x)}  \label{4.4}
\end{eqnarray}
for $f \in L^1(I, \mu)$, where $\hat{f}(x):=(x+N)f(x)h(x)$, $x \in I$.
In addition, the $n$th power $S^n$ of \ $S$ is written as follows:
\begin{equation}
S^n f(x) = \frac{U^n \hat{f}(x)}{(x+N)h(x)} \label{4.5}
\end{equation}
for any $f \in L^1(I, \mu)$ and any $n \geq 1$.
  \item[(b)]
Let $K$ denote the Perron-Frobenius operator of \ $T_N$ under $\lambda$. Then the following holds a.e. in $I$:
\begin{equation}
K f(x) = \sum_{i \geq N} \frac{N}{(x+i)^2} f\left(\frac{N}{x+i}\right), \ f \in L^1(I, \lambda). \label{4.6}
\end{equation}
In addition, the $n$th power $K^n$ of \ $K$ is written as follows:
\begin{equation}
K^n f(x) = \frac{U^n \hat{f}(x)}{x+N}, \ f \in L^1(I, \lambda), \label{4.7}
\end{equation}
for any $f \in L^1(I, \lambda)$ and any $n \geq 1$, where $\hat{f}(x):=(x+N)f(x)$, $x \in I$.
  \item[(c)]
For any $n \in \mathbb{N}_+$ and $A \in {\mathcal{B}}_{I}$,
we have
\begin{equation}
\mu \left((T_{N})^{-n}(A)\right)
= \int_{A} U^nf(x) d G_{N}(x) \label{4.8}
\end{equation}
where $f(x):= (\log\left(\frac{N+1}{N}\right)) (x+N) h(x)$ for
$x \in I$.
\end{enumerate}
\end{enumerate}
\end{proposition}

\noindent \textbf{Proof.}
(i) Let $T_{N,i}$ denote the restriction of $T_N$ to the subinterval
$I_i:=\left(\frac{N}{i+1}, \frac{N}{i}\right]$, $i \geq N$, that is,
\begin{equation}
T_{N,i}(x) = \frac{N}{x} - 1, \quad x \in I_i. \label{4.9}
\end{equation}
Let $C(A):=\left(T_{N}\right)(A)$ and $C_{i}(A):=\left(T_{N,i}\right)^{-1}(A)$ for
$A \in {\mathcal B}_I$.
Since $C(A)=\bigcup_{i}C_i(A)$ and $C_i\cap C_j$ is a null set when $i \neq j$,
we have
\begin{equation}
\int_{C(A)} f \,d G_N = \sum_{i \geq N} \int_{C_i(A)}f\, d G_N,
\quad
f \in L^1(I,G_N),\,A \in {\mathcal{B}}_I. \label{4.10}
\end{equation}
For any $i \geq N$, by the change of variables
$x = \left(T_{N,i}\right)^{-1}(y) = \displaystyle \frac{N}{y+i}$,
we successively obtain
\begin{eqnarray}
\int_{C_i(A)}f(x) \,G_N(dx) &=& \left(\log\left(\frac{N+1}{N}\right)\right)^{-1} \int_{C_i(A)} \frac{f(x)}{N+x}\,dx \nonumber \\
\nonumber \\
&=& \left(\log\left(\frac{N+1}{N}\right)\right)^{-1} \int_{A}\frac{f\left(\frac{N}{y+i}\right)}{N+\frac{N}{y+i}}
\frac{N}{(y+i)^2}dy \nonumber \\
\nonumber \\
&=& \int_{A} V_{N,i}(y)\, f\left(\frac{N}{y+i}\right)\,G_N (dy). \label{4.11}
\end{eqnarray}
Now, (\ref{4.2}) follows from (\ref{4.10}) and (\ref{4.11}).

\noindent
(ii)(a) From (\ref{4.9}), for any $f \in L^1(I,G_N)$
and  $A \in {\mathcal{B}}_I$, we have
\begin{equation} \label{4.12}
\displaystyle \int_{C(A)} f(x) \,\mu(dx)
= \displaystyle \sum_{i \geq N} \displaystyle \int_{C_i(A)}f(x) \,\mu(dx).
\end{equation}
Then
\begin{equation} \label{4.13}
\int_{C_i(A)}f(x) \,\mu(dx)
= \displaystyle \displaystyle \int_{C_i(A)}f(x)h(x)\, dx
= \displaystyle \displaystyle \int_{A} \displaystyle f\left(\frac{N}{y+i}\right)\,h\left(\frac{N}{y+i}\right)\frac{N}{(y+i)^2} \,dy.
\end{equation}
From (\ref{4.12}) and (\ref{4.13}),
\begin{equation} \label{4.14}
\displaystyle
\displaystyle \int_{C(A)} f(x) \,\mu(dx)
= \int_{A} \sum_{i \geq N} f\left(\frac{N}{x+i}\right)\,h\left(\frac{N}{x+i}\right)\frac{N}{(x+i)^2} \,dx.
\end{equation}
Since $d\mu = h dx$, (\ref{4.3}) follows from (\ref{4.14}). Now, since $\hat{f}(x)=(x+N)f(x)h(x)$, from (\ref{4.3}), we have
\begin{equation} \label{4.15}
U \hat{f}(x) = N(x+N)
\sum_{i \geq N}\frac{h\left(\frac{N}{x+i}\right)}{(x+i)^2}f\left(\frac{N}{x+i}\right).
\end{equation}
From (\ref{4.3}) and (\ref{4.15}), (\ref{4.4}) follows immediately.

\noindent
(ii)(b)\, The formula (\ref{4.6}) is a consequence of (\ref{4.4}) and follows immediately.

\noindent
(ii)(c)\, We will use mathematical induction.
For $n=0$, the equation (\ref{4.8}) holds by definitions of $f$ and $h$.
Assume that (\ref{4.8}) holds for some $n \in \mathbb{N}$. Then
\begin{equation} \label{4.16}
\mu \left(\left(T_N\right)^{-(n+1)}(A)\right) =
\mu \left(\left(T_N\right)^{-n}\left(\left(T_N\right)^{-1}(A)\right)\right)
= \int_{C(A)} U^n f(x)\,G_N(dx).
\end{equation}
Since $U=U_{T_N}$ and (\ref{4.1}), we have
\begin{equation} \label{4.17}
\int_{C(A)} U^n f(x) \,G_N(dx) = \int_{A} U^{n+1} f(x) \,G_N(dx).
\end{equation}
Therefore,
\begin{equation} \label{4.18}
\mu \left(\left(T_N\right)^{-(n+1)}(A)\right)
= \int_{A} U^{n+1} f(x)G_N(dx)
\end{equation}
which ends the proof.

\hfill $\Box$

For a function $f:I \to {\mathbb C}$, define the {\it variation} ${\rm var}_{A}f$
of $f$ on a subset $A$ of $I$ by
\begin{equation}
{\rm var}_A f := \sup \sum^{k-1}_{i=1} |f(t_{i+1}) - f(t_{i})| \label{4.19}
\end{equation}
where the supremum being taken over $t_1 < \cdots < t_k$, $t_i \in A$,
$1 \leq i \leq k$, and $k \geq 2$ (\cite{IK-2002}, p75).
We write simply $\mathrm{var} f$
for $\mathrm{var}_I f$.
Let $BV(I):=\{f:I\to {\mathbb C}:
{\rm var}\,f<\infty\}$ and let $L^{\infty}(I)$ denote
the collection of all bounded measurable functions $f:I \rightarrow \mathbb{C}$.
It is known that $BV(I) \subset L^{\infty}(I) \subset L^1(I,\mu)$.
Let  $L(I)$ denote the Banach space of all complex-valued Lipschitz continuous functions on $I$
with the following norm $\|\cdot\|_L$:
\begin{equation}
\left\| f \right\|_L := \sup_{x \in I} |f(x)| + s(f), \label{4.20}
\end{equation}
with
\begin{equation}
s(f):=\sup_{x \ne y} \frac{|f(x) - f(y)|}{|x - y|}, \quad f \in L(I). \label{4.21}
\end{equation}

In the following proposition we show that the operator $U$ in (\ref{4.2}) preserves monotonicity and enjoys a contraction property for Lipschitz continuous functions.
\begin{proposition} \label{prop4.2}
Let $U$ be as in $(\ref{4.2})$.
\begin{enumerate}
  \item[(i)]
  Let $f \in L^{\infty}(I)$. Then the following holds:
\begin{enumerate}
    \item[(a)]
	If $f$ is non-decreasing (non-increasing),
	then $Uf$ is non-increasing (non-decreasing).
    \item[(b)]
	If $f$ is monotone, then
\begin{equation}
\mathrm{var}\,(Uf) \leq \frac{1}{N+1} \cdot \mathrm{var } f. \label{4.22}
\end{equation}
\end{enumerate}

 \item[(ii)]
 For any $f \in L(I)$, we have
\begin{equation}
s(Uf) \leq q \cdot s(f), \label{4.23}
\end{equation}
where
\begin{equation}
q:=N \left( \sum_{i \geq N} \left(\frac{N}{i^3(i+1)}+ \frac{i+1-N}{i(i+1)^3}\right) \right). \label{4.24}
\end{equation}
\end{enumerate}
\end{proposition}

\noindent \textbf{Proof.}
(i)(a)
To make a choice assume that $f$ is non-decreasing.
Let $x < y$, \, $x, y \in I$.
We have $U f(y) - U f(x) = S_1 + S_2$, where
\begin{eqnarray}
S_1 &=& \sum_{i \geq N} V_{N,i}(y) \left(f\left(\frac{N}{y+i}\right) - f\left(\frac{N}{x+i}\right)\right), \label{4.25} \\
S_2 &=& \sum_{i \geq m} \left(V_{N,i}(y) - V_{N,i}(x)\right)f\left(\frac{N}{x+i}\right). \label{4.26}
\end{eqnarray}
Clearly, $S_1 \leq 0$. Now, since $\sum_{i \geq N} V_{N,i}(x) = 1$ for any $x \in I$,
we can write
\begin{equation}
S_2 = - \sum_{i \geq N} \left( f\left(\frac{N}{x+N}\right) - f\left(\frac{N}{x+i}\right) \right) \left(V_{N,i}(y) - V_{N,i}(x)\right). \label{4.27}
\end{equation}
As is easy to see, the functions $V_{N,i}$ are increasing for all $i \geq N$.
Also, using that $f\left(\frac{N}{x+N}\right) \geq f\left(\frac{N}{x+i}\right)$, we have that $S_2 \leq 0$. Thus $U f(y) - U f(x) \leq 0$ and the proof is complete.

\noindent
(i)(b)
Assume that $f$ is non-decreasing. Then by (a) we have
\begin{equation}
\mathrm{var}\, Uf = Uf(0) - Uf(1) = \sum_{i \geq N} \left( V_{N,i}(0) f\left(\frac{N}{i}\right) - V_{N,i}(1) f\left(\frac{N}{1+i}\right) \right). \label{4.28}
\end{equation}
By calculus, we have
\begin{eqnarray}
\mathrm{var}\, Uf &=& \sum_{i \geq N}
\left( \frac{N}{i(i+1)}f\left(\frac{N}{i}\right) - \frac{N+1}{(i+1)(i+2)}f\left(\frac{N}{i+1}\right)\right) \nonumber\\
&=& \frac{1}{N+1} f(1) - \sum_{i \geq N} \frac{1}{(i+1)(i+2)}f\left(\frac{N}{i+1}\right) \nonumber \\
&\leq& \frac{1}{N+1} f(1) - \sum_{i \geq N} \left(\frac{1}{i+1}-\frac{1}{i+2}\right)f(0)  \nonumber \\
&=& \frac{1}{N+1} (f(1) - f(0)) = \frac{1}{N+1}\mathrm{var} f. \nonumber
\end{eqnarray}
\noindent
(ii) For $x \ne y$, $x,y \in I$, we have
\begin{eqnarray}
\frac{Uf(y) - Uf(x)}{y - x} &=& \sum_{i \geq N} \frac{V_{N,i}(y) - V_{N,i}(x)}{y - x}f\left(\frac{N}{x+i}\right) \nonumber \\
&-& \sum_{i \geq N} V_{N,i}(y) \frac{f\left(\frac{N}{y+i}\right) - f\left(\frac{N}{x+i}\right)}{\frac{N}{x+i} - \frac{N}{x+i}}\cdot \left(\frac{N}{x+i}\right)\left(\frac{N}{y+i}\right). \label{4.29}
\end{eqnarray}
Remark that
\begin{equation}
V_{N,i}(u) = \frac{i+1-N}{u+i+1} + \frac{N-i}{u+i}, \quad i \geq N, \label{4.30}
\end{equation}
and then
\begin{eqnarray}
&{}& \sum_{i \geq N}
\frac{V_{N,i}(y) - V_{N,i}(x)}{y - x}f\left(\frac{N}{x+i}\right) \nonumber \\
&=& \sum_{i \geq N} \frac{i+1-N}{(y + i+1)(x + i+1)} \left( f\left(\frac{N}{x+i+1}\right) - f\left(\frac{N}{x+i}\right) \right). \label{4.31}
\end{eqnarray}

Assume that $x > y$. It then follows from (\ref{4.31}) and (\ref{4.29}) that
\begin{eqnarray}
\left|\frac{Uf(y) - Uf(x)}{y - x}\right| &\leq& s(f) \sum_{i \geq N}
\left( \frac{N(i+1-N)}{(y+i)(y + i+1)^3} + \frac{N \cdot V_{N,i}(y)}{(y+i)^2}\right) \nonumber \\
&\leq& q \cdot s(f) \label{4.32}
\end{eqnarray}
where $q$ is as in (\ref{4.24}).
Since
\begin{equation}
s(Uf) = \sup_{x,y \in I, x \geq y} \left|\frac{Uf(y) - Uf(x)}{y - x}\right| \label{4.33}
\end{equation}
then the proof is complete.
\hfill $\Box$
\\

\noindent
\textbf{Proof of Proposition \ref{prop2.9}}
For $(I,{\mathcal B}_I,G_N,T_N)$ in Definition \ref{def2.1}(ii),
let $U$ denote its Perron-Frobenius operator.
Let
\begin{equation}
\rho_N(x):=\frac{k_N}{x+N}, \quad x \in I, \label{4.34}
\end{equation}
where $k_N=\left(\log\left(\frac{N+1}{N}\right)\right)^{-1}$.
From properties of the Perron-Frobenius operator,
it is sufficient to show that the function $\rho_N$ defined in (\ref{4.34})
satisfies $U\rho_N=\rho_N$.

Since $\left(T_N\right)^{-1}(x)=\left\{\displaystyle \frac{N}{x+i}, i \geq N, x \in I \right\}$, we have
\begin{eqnarray}
  U \rho_N(x) &=& \frac{d}{dx}\int_{\left(T_N\right)^{-1}([0,x])}\rho_N(t) \,dt = \sum_{t \in \left(T_N\right)^{-1}(x)}\frac{\rho_N(t)}{\left|(T_N)'(t)\right|}\nonumber \\
              &=&\sum_{i \geq N} \frac{N}{(x+i)^2}\, \rho_N\left(\frac{N}{x+i}\right). \label{4.35}
\end{eqnarray}

By definition of $\rho_N$, we see that
\begin{equation} \label{4.36}
U \rho_N(x) = \sum_{i \geq N}\frac{1}{(x+i)(x+i+1)}
= \rho_N(x).
\end{equation}
Hence the statement is proved.
\hfill $\Box$

\end{document}